\documentclass[
11pt,twoside, a4paper,
]{amsart}

\usepackage[all]{xy}
\usepackage{mathptmx}
\usepackage{amscd}
\usepackage{amssymb}
\usepackage{eucal}
\usepackage{amsxtra}
\usepackage{verbatim}
\usepackage{url}

\usepackage[english,francais]{babel}

\title[Conjugate homogeneous spaces]
{           {\protect\hfill \normalfont \tiny
            \\ \vspace{10pt}}
Conjugate complex  homogeneous spaces \\
with non-isomorphic fundamental groups 
}

\author{Mikhail Borovoi}

\address{Borovoi:
Raymond and Beverly Sackler School of Mathematical Sciences,
Tel Aviv University, 6997801 Tel Aviv, Israel}
\email{borovoi@post.tau.ac.il}

\author{Yves  Cornulier}

\address{Cornulier:
D\'epartement de Math\'ematiques, B\^atiment 425, Universit\'e Paris-Sud 11, 91405 Orsay, France}
\email{yves.cornulier@math.u-psud.fr}

\thanks{M. B. was partially supported by the Hermann Minkowski Center for Geometry}
\thanks{Y. C. was supported by ANR GSG 12-BS01-0003-01}

\date{May 9, 2015}


\usepackage{comment}

\theoremstyle{plain}

\newtheorem{theorem}{Theorem} [section]
\newtheorem{proposition}[theorem]{Proposition}
\newtheorem{lemma}[theorem]{Lemma}
\newtheorem{corollary}[theorem]{Corollary}

\newtheorem{conditional-result}[theorem]{Conditional Result}

\newtheorem{theorem?}{Theorem(?)} [section]
\newtheorem{proposition?}[theorem]{Proposition(?)}
\newtheorem{lemma?}[theorem]{Lemma(?)}
\newtheorem{corollary?}[theorem]{Corollary(?)}

\newtheorem*{theorem*}{Theorem}
\newtheorem*{proposition*}{Proposition}
\newtheorem*{lemma*}{Lemma}
\newtheorem*{corollary*}{Corollary}
\newtheorem*{question*}{Question}
\newtheorem*{conjecture*}{Conjecture}
\newtheorem*{claim*}{Claim}

\newtheorem*{introtheorem*}{Theorem}
\newtheorem*{introproposition*}{Proposition}
\newtheorem*{introlemma*}{Lemma}
\newtheorem*{introcorollary*}{Corollary}

\theoremstyle{definition}

\newtheorem*{definition*}{Definition}
\newtheorem*{example*}{Example}

\theoremstyle{remark}

\newtheorem*{remark*}{Remarque}

\DeclareSymbolFont{rsfs}{U}{rsfs}{m}{n}
\DeclareSymbolFontAlphabet{\mathcal}{rsfs}

\DeclareFontEncoding{OT2}{}{} 
\DeclareTextFontCommand{\textcyr}{\fontencoding{OT2}
    \fontfamily{wncyr}\fontseries{m}\fontshape{n}\selectfont}

\newcommand{\isoto}{\overset{\sim}{\longrightarrow}}

\newcommand{\into}{\hookrightarrow}

\newcommand{\labelto}[1]{\xrightarrow{\makebox[1.5em]{\scriptsize ${#1}$}}}

\newcommand{\C}{{\mathbf{C}}}
\newcommand{\Z}{{\mathbf{Z}}}
\newcommand{\Q}{{\mathbf{Q}}}

\renewcommand{\ker}{{\rm ker}}

\newcommand{\im}{{\rm im}}

\newcommand{\Aut}{{\rm Aut}}
\def\Out{{\rm Out}}

\def\top{{\textup{top}}}

\def\SL{{\rm SL}}

\def\tiltheta{{\tilde{\theta}}}

\def\vk{{\kappa}}
\def\vkbar{{\overline{\vk}}}

\def\ve{\varepsilon}
\def\rth{\tfrac{1}{r}}

\def\hs{\kern 0.6pt}

\begin{document}

\selectlanguage{english}
\date{\today}

\begin{abstract}
Let $X=G/\Gamma$ be the quotient of a connected reductive algebraic $\C$-group $G$  by a finite subgroup $\Gamma$.
We describe the topological fundamental group of the homogeneous space $X$, which is nonabelian when $\Gamma$ is nonabelian.
Further, we construct an example of a homogeneous space $X$ and an automorphism $\sigma$ of $\C$
such that the topological fundamental groups of $X$ and of the conjugate variety $\sigma X$
are not isomorphic.
\bigskip

\selectlanguage{francais}
\noindent{\sc R\'esum\'e.} \  {\em Espaces homog\`enes complexes conjugu\'es avec  groupes fondamentaux non isomorphes.} \
Soit $X=G/\Gamma$ le quotient d'un $\C$-groupe alg\'ebrique r\'eductif connexe $G$  par un sous-groupe fini $\Gamma$.
On d\'ecrit le groupe fondamental topologique de l'espace homog\`ene $X$, qui est non ab\'elien quand $\Gamma$ est non ab\'elien.
Puis on construit un exemple d'espace homog\`ene $X$ et d'automorphisme $\sigma$ de $\C$ tels que les groupes fondamentaux
topologiques de $X$ et de la vari\'et\'e conjugu\'ee $\sigma X$ ne sont pas isomorphes.
\end{abstract}

\maketitle

\selectlanguage{francais}
\section*{Abridged French version}
Soit $X$ une vari\'et\'e alg\'ebrique point\'ee d\'efinie sur le corps $\C$ des nombres complexes, suppos\'ee irr\'eductible et quasi-projective.
L'espace topologique point\'e $X(\C)$ est alors connexe;
on d\'esigne par $\pi_1(X):=\pi_1^\top(X(\C))$ son groupe fondamental, appel\'e groupe fondamental topologique de $X$.
Soit  $\sigma$  un automorphisme du corps $\C$ (pas forc\'ement continu).
En appliquant $\sigma$ aux coefficients des polyn\^omes d\'efinissant $X$, on obtient une vari\'et\'e $\sigma X$ sur $\C$, dite  vari\'et\'e conjugu\'ee.
Les compl\'et\'es profinis des groupes $\pi_1(X)$ et $\pi_1(\sigma X)$ sont canoniquement isomorphes (comme groupes topologiques),
car ils s'identifient naturellement au groupe fondamental \'etale de $X$.
En revanche, les groupes $\pi_1(X)$ et $\pi_1(\sigma X)$ ne sont pas toujours isomorphes,
par un r\'esultat de Serre \cite{Serre}.
Les exemples de Serre comprennent des surfaces projectives lisses.
D'autres exemples ont \'et\'e obtenus plus r\'ecemment: des vari\'et\'es de Shimura dans \cite{MS, R},
et des surfaces projectives dans \cite{BCG, GJ} pour des choix tr\`es g\'en\'eraux de l'automorphisme $\sigma$
(dans \cite{GJ} pour tout $\sigma$
dont la restriction \`a $\overline{\Q}$  diff\`ere de l'identit\'e et de la conjugaison complexe).

Dans cette note, nous donnons un exemple des {\em espaces homog\`enes} conjugu\'es avec groupes fondamentaux topologiques non isomorphes.
Le plan de la note est le suivant.
Nous consid\'erons, dans le \S\ref{s:fundam}, les groupes fondamentaux de certains espaces homog\`enes topologiques de la forme $G/\Gamma$,
o\`u $G$ est un groupe de Lie r\'eel connexe et $\Gamma\subset G$ est un sous-groupe discret.
Nous en d\'eduisons, dans le \S\ref{s:algebraic}, une formule explicite pour d\'ecrire le groupe fondamental $\pi_1(G/\Gamma)$ dans le cas o\`u
$G$ est un groupe alg\'ebrique lin\'eaire connexe d\'efini sur $\C$, et $\Gamma$  est un sous-groupe fini de $G$.
En utilisant cette formule, nous construisons dans le \S\ref{s:example}
un exemple d'espace homog\`ene affine $X=G/\Gamma$  d\'efini sur $\C$
et un automorphisme $\sigma$ de $\C$
tels que les groupes fondamentaux topologiques $\pi_1(\sigma X)$ et $\pi_1(X)$ ne sont pas isomorphes.
Pr\'ecis\'ement, on choisit $G=\SL(n,\C)\times\C^*$ avec $n\ge 5$, et $\Gamma$  un sous-groupe non ab\'elien fini d'ordre 55.
L'inclusion de $\Gamma$ dans $G$ est donn\'ee
par un plongement arbitraire de $\Gamma$ dans $\SL(n,\C)$ et  par un homomorphisme non trivial de $\Gamma$ dans $\C^*$.
Notre formule permet de v\'erifier que $\pi_1(X)$ et isomorphe \`a $(\Z/11\Z)\rtimes_4\Z$,
o\`u la notation signifie que le g\'en\'erateur $1$ de $\Z$ agit sur $\Z/11\Z$ par multiplication par 4,
tandis que pour  $\sigma$ envoyant $\zeta=\exp 2\pi i/5$ sur $\zeta^2$, le groupe fondamental $\pi_1(\sigma X)$
de la vari\'et\'e conjugu\'ee est isomorphe \`a $(\Z/11\Z)\rtimes_9\Z$.
Un argument simple permet de v\'erifier que ces deux groupes ne sont pas isomorphes.

\selectlanguage{english}

\section{Introduction}
\label{s:Intro}

Let $X$ be a pointed algebraic variety defined over  $\C$.
We assume that $X$ is irreducible and quasi-projective.
The pointed topological space $X(\C)$ is then connected, and
we denote by $\pi_1(X)$ the topological fundamental group of $X(\C)$, i.e.,
$\pi_1(X):=\pi_1^\top(X(\C))$.
Let $\sigma$ be a field automorphism of $\C$, not necessarily continuous.
On applying $\sigma$ to the coefficients of the polynomials defining
$X$, we obtain a conjugate algebraic variety $\sigma X$ over $\C$.
Though the profinite completions of $\pi_1(X)$ and $\pi_1(\sigma X)$ are isomorphic,
the groups $\pi_1(X)$ and $\pi_1(\sigma X)$ themselves are not necessarily isomorphic.
Serre \cite{Serre} obtained the first examples of conjugate varieties $X$ and $\sigma X$
with $\pi_1(\sigma X)\not\simeq\pi_1(X)$.
Serre's examples include smooth projective surfaces.
More examples were obtained recently: Shimura varieties in \cite{MS} and \cite{R},
and smooth projective surfaces in \cite{BCG} and \cite{GJ} for a very general choice of $\sigma$
(in \cite{GJ} for any $\sigma$ whose restriction to $\overline{\Q}$ differs from the identity and the complex conjugation).

In this note we give  an example of conjugate  {\em homogeneous spaces} with non-isomorphic topological fundamental groups.
The outline of the note is as follows.
In Section \ref{s:fundam} we consider topological homogenous spaces of the form $G/\Gamma$,
where  $G$ is a connected real Lie group and $\Gamma\subset G$ is a discrete subgroup.
In Section \ref{s:algebraic} we write an explicit formula for $\pi_1(G/\Gamma)$
when $G$ is a complex linear algebraic group and $\Gamma\subset G$ is a finite subgroup.
In Section \ref{s:example} using this formula we construct an example
of an affine homogeneous space $X=G/\Gamma$ over $\C$ and an automorphism $\sigma$ of $\C$
such that $\pi_1(\sigma X)$ is not isomorphic to $\pi_1(X)$.
In our example $G=\SL(n,\C)\times \C^*$ with $n\ge 5$, and $\Gamma$ is a nonabelian finite subgroup of order 55.

\section{The quotient of a Lie group by a discrete subgroup}
\label{s:fundam}

Let
\[1\to S\labelto{i} G\labelto{\tau} T\to 1\]
be a short exact sequence of connected real Lie groups.
Let $\Gamma\subset G$ be a discrete subgroup such that the projection $\Lambda=\tau(\Gamma)\subset T$ is discrete.
Our goal is to describe $\pi_1(G/\Gamma)$, where $G/\Gamma$ is viewed as a pointed manifold with base point the image of 1.

Set  $\Gamma_S=\Gamma\cap S$.
The homomorphism $\tau\colon G\to T$ induces a fibration $ G/\Gamma\to T/\Lambda$ with fiber $S/\Gamma_S$\hs,
which gives rise to an  exact sequence in homotopy groups
\[\pi_1(S/\Gamma_S)\labelto{i_{*}}\pi_1(G/\Gamma)\labelto{\tau_{*}} \pi_1(T/\Lambda)\to 1.\]
The fibration $G\to G/\Gamma$ with fiber $\Gamma$ gives rise to an exact sequence in homotopy groups
\[1\to\pi_1(G)\to\pi_1(G/\Gamma)\labelto{f} \Gamma\to 1,\]
where $f$ is a homomorphism by Lemma \ref{l:connecting} below.
Considering the above  fibrations and also the fibrations $S\to S/\Gamma_S$, $T\to T/\Lambda$ and $G\to T$,
we obtain the following commutative diagram of groups and homomorphisms with exact rows and columns:
\[
\xymatrix@R=18pt{
1\ar[r] &\pi_1(S)\ar[d]\ar[r] &\pi_1(S/\Gamma_S)\ar[d]^(0.44){i_{*}}\ar[r]    &\Gamma_S\ar[d]^(0.44){i} \ar[r] &1 \\
1\ar[r] &\pi_1(G)\ar[d]\ar[r] &\pi_1(G/\Gamma)\ar[d]^(0.44){\tau_{*}}\ar[r]^-f &\Gamma\ar[d]^(0.44){\tau}\ar[r] &1 \\
1\ar[r] &\pi_1(T)\ar[d]\ar[r] &\pi_1(T/\Lambda)\ar[d]\ar[r]^-{f_T}       &\Lambda\ar[d]\ar[r]     &1\\
        &1                    &1                                        &1
}
\]
From this diagram we obtain homomorphisms
\[\chi\colon \pi_1(S)\to\pi_1(S/\Gamma_S)\labelto{i_*} \pi_1(G/\Gamma)\quad\text{and}\quad
\phi\colon\pi_1(G/\Gamma)\to \pi_1(T/\Lambda)\underset{\Lambda}{\times} \Gamma,\]
where the fiber product $\pi_1(T/\Lambda)\times_\Lambda\hs \Gamma$ is the group of pairs $(x,\gamma)\in \pi_1(T/\Lambda)\times \Gamma$
such that $f_T(x)=\tau(\gamma)$.
The homomorphism $\phi$ takes $y\in\pi_1(G/\Gamma)$ to the pair $(\tau_*(y),f(y))\in \pi_1(T/\Lambda)\times_\Lambda\hs \Gamma$.

\begin{theorem}\label{t:pi-1}
With the above notation, the sequence
\[\pi_1(S)\labelto{\chi}\pi_1(G/\Gamma)\labelto{\phi} \pi_1(T/\Lambda)\underset{\Lambda}\times \Gamma\to 1\]
is exact. In particular, if $S$ is simply connected, then $\phi$ is an isomorphism.
\end{theorem}

\noindent {\bf Proof.}
We prove the theorem by diagram chasing.
Clearly $\phi\circ\chi=1$. We show that $\ker\,\phi\subset\im\,\chi$.
Let $y\in\ker\,\phi\subset \pi_1(G/\Gamma)$, then $f(y)=1$ and  $\tau_*(y)=1$.
Then $y$ comes from some element $z\in\pi_1(G)$, whose image in $\pi_1(T)$ is 1.
Hence $z$ comes from some element $u\in\pi_1(S)$.
We see that $y=\chi(u)$, as required.

We show that $\phi$ is surjective.
Let $(x,\gamma)\in  \pi_1(T/\Lambda)\times_\Lambda \Gamma$, i.e.,
$x\in\pi_1(T/\Lambda)$, $\gamma\in\Gamma$, and $f_T(x)=\tau(\gamma)$.
We can lift $x$ to some  element $y\in\pi_1(G/\Gamma)$,
then $\tau(f(y))=\tau(\gamma)$.
Set $z=f(y)\gamma^{-1}$, then $\tau(z)=1$, hence $z$ comes from some element of $\Gamma_S$
and from some element $u$ of $\pi_1(S/\Gamma_S)$.
Set $y'=i_*(u)^{-1}y\in\pi_1(G/\Gamma)$, then $f(y')=\gamma$ and $\tau_*(y')=\tau_*(y)=x$.
We see that $(x,\gamma)=\phi(y')$, as required.
\qed
\medskip

The following lemma, which we used above, is well-known.

\begin{lemma}\label{l:connecting}
Let $G$ be a connected Lie group, $\Gamma\subset G$ be a (closed) Lie subgroup, not necessarily connected.
Then the connecting map $f\colon \pi_1(G/\Gamma)\to \pi_0(\Gamma)$ in the exact sequence
\[\pi_1(\Gamma)\to\pi_1(G)\to\pi_1(G/\Gamma)\labelto{f}\pi_0(\Gamma)\to 1\]
is a homomorphism.
\end{lemma}

\noindent {\bf Proof.}
Denote by $\lambda\colon \Gamma\to\pi_0(\Gamma)$ the canonical epimorphism.
Consider two based loops $\theta_i\colon[0,1]\to G/\Gamma$ in $G/\Gamma$  ($i=1,2$).
Let $\tiltheta_i\colon [0,1]\to G$ be a path lifting the loop $\theta_i$ to $G$
with $\tiltheta_i(0)=1$, and set $\gamma_i=\tiltheta_i(1)\in \Gamma$.
By definition $f([\theta_i])=\lambda(\gamma_i)\in\pi_0(\Gamma)$, where $[\theta_i]$
denotes the class of the based loop $\theta_i$ in $\pi_1(G/\Gamma)$.
Then $\gamma_1\tiltheta_2$ is a path in $G$ from $\gamma_1$ to $\gamma_1\gamma_2$ mapping in $G/\Gamma$ to the loop $\theta_2$,
hence the concatenation of $\tiltheta_1$ and $\gamma_1\tiltheta_2$ is a path in $G$ from $1$ to $\gamma_1\gamma_2$
mapping in $G/\Gamma$ to the loop obtained by concatenation of $\theta_1$ and $\theta_2$.
Thus $f([\theta_1]\cdot[\theta_2])=\lambda(\gamma_1 \gamma_2)=\lambda(\gamma_1)\,\lambda(\gamma_2)=f([\theta_1])\, f([\theta_2])$, as required.
\qed

\section{The quotient of a complex algebraic group by a finite subgroup}
\label{s:algebraic}

Let $G$ be a connected linear algebraic group over $\C$.
Let $\Gamma\subset G$ be a finite subgroup.
Set $X=G/\Gamma$.
We wish to compute the topological fundamental group $\pi_1(X)$.

Let $U$ denote the unipotent radical of $G$, then $G':=G/U$ is reductive.
The canonical epimorphism $\rho\colon G\to G'$
induces a fibration $G/\Gamma\to G'/\Gamma'$ with fiber $U$, where $\Gamma'=\rho(\Gamma)$, and hence,
the induced homomorphism $\rho_*\colon\pi_1(G/\Gamma)\to\pi_1(G'/\Gamma')$ is an isomorphism.
Therefore, we may and shall assume that $G$ is reductive.
Replacing the reductive group $G$ by a finite cover and $\Gamma$ by its inverse image, we may and shall
assume that the semisimple group $S:=[G,G]$ is simply connected.
Let $\Lambda$ denote the image of $\Gamma$ in the algebraic torus  $T:=G/S$,
then $T/\Lambda$ is also an algebraic torus, hence $\pi_1(T/\Lambda)$ is a free abelian group isomorphic to $\Z^{{\rm dim}\, T}$.
The next corollary, which follows immediately from Theorem \ref{t:pi-1},
describes $\pi_1(G/\Gamma)$ in terms of $\Gamma$ and the free abelian group $\pi_1(T/\Lambda)$.

\begin{corollary}\label{c:pi-1-alg}
Let $G$ be a connected reductive algebraic group over $\C$ such that the commutator subgroup $S$ of $G$ is simply connected.
Set $T=G/S$.
Let $\Gamma\subset G$ be a finite subgroup, and let $\Lambda$ denote the image of $\Gamma$ in $T$.
Then there is a canonical isomorphism
\[ \pi_1(G/\Gamma)\isoto \pi_1(T/\Lambda)\underset{\Lambda}\times \Gamma,\]
where $\pi_1(T/\Lambda)\times_{\Lambda}\hs \Gamma$ is the fiber product
with respect to the epimorphism $\pi_1(T/\Lambda)\to\Lambda$ of Lemma \ref{l:connecting}
and the canonical epimorphism $\Gamma\to\Lambda$.
\end{corollary}

\section{Example}
\label{s:example}

Let $A=\Z/m\Z$, the additive group of residues modulo $m$.
Let $B\subset (\Z/m\Z)^*$ be a {\em cyclic} subgroup of some order $r$
in the multiplicative group of invertible residues modulo $m$.
The group $B$ acts naturally on $A$ by multiplication:
an element $b\in B\subset(\Z/m\Z)^*$ acts by  $a\mapsto ba$.
Set
\[H=A\rtimes B \]
(the semidirect product).
We regard $B$ as a subgroup of $H$.
Consider an embedding $\varphi\colon B\into\C^*$, then $\varphi(B)=\mu_r\subset \C^*$, the group of $r$-th roots of unity.

Choose an embedding $\alpha\colon H\into \SL(n,\C)$ for some natural number $n$. Set
\[G= \SL(n,\C)\times \C^*.\]
For $(a,b)\in A\rtimes B=H$ set
\[\psi(a,b)=(\alpha(a,b),\varphi(b))\in \SL(n,\C)\times\C^*.\]
We obtain an embedding $\psi=\psi_{\alpha,\varphi}\colon H\into G$.
Set $\Gamma=\psi(H),$ \ $X=X_{\alpha,\varphi}=G/\Gamma$.
Then $X$ is an affine algebraic variety over $\C$.

Let $b\in B$.
Write $A\rtimes_b \Z$ for the semidirect product of $A$ and $\Z$,
where the generator $1$ of $\Z$ acts on $A$ by multiplication by $b$.
Set $\zeta=\exp 2\pi i/r\in \mu_r$.

\begin{proposition}\label{p:pi-1-X}
$\pi_1(X_{\alpha,\varphi})\simeq (\Z/m\Z)\rtimes_{\varphi^{-1}(\zeta)}\Z$.
\end{proposition}

\noindent {\bf Proof.}
Set $S=\SL(m,\C)$, $T=\C^*$.
Let $\tau\colon G=S\times T\to T$ denote the projection, then
$\tau(\psi(a,b))=\varphi(b)$ for $(a,b)\in A\rtimes B =H$.
Set $\Lambda=\tau(\Gamma)=\tau(\psi(H))\subset T$,
then $\Lambda=\varphi(B)=\mu_r\subset\C^*=T$.

Consider the following universal covering of $T=\C^*$:
\[ \ve\colon \C\to\C^*=T,\quad z\mapsto \exp 2\pi iz \ \text{ for }\ z\in\C,\]
it induces a universal covering of $T/\Lambda$:
\[ \C\labelto{\ve}\C^*\to\C^*/\mu_r=T/\Lambda\simeq\C^*.\]
We identify $\pi_1(T/\Lambda)$ with $\ve^{-1}(\mu_r)=\rth\Z\subset\C$,
then the homomorphism $\pi_1(T/\Lambda)\to\Lambda=\mu_r$ of Lemma \ref{l:connecting}
is the restriction of $\ve$ to $\rth\Z$, hence it takes
the generator $\rth\in\rth\Z=\pi_1(T/\Lambda)$ to the element $\ve(\rth)=\zeta\in\mu_r$.

Since $S=\SL(n,\C)$ is simply connected, by Corollary \ref{c:pi-1-alg} we have
\[\pi_1(X_{\alpha,\varphi})=\pi_1(G/\Gamma)=\pi_1(T/\Lambda)\underset{\Lambda}{\times}\Gamma\simeq\rth\Z\underset{\mu_r}{\times}H,\]
where the homomorphism $\rth\Z\to\mu_r$ takes $\rth$ to $\zeta$
and the homomorphism $H\to\mu_r$ takes $(a,b)\in H$ to $\tau(\psi(a,b))=\varphi(b)$.
Since $\rth\Z$ is a free abelian group, the group extension
\[ 1\to \{0\}\times A\to \rth\Z\underset{\mu_r}{\times}H\to \rth\Z\to 1\]
splits, hence $\pi_1(X_{\alpha,\varphi})\simeq A\rtimes\rth\Z$.
The action of $\rth\Z$ on $A$ in this semidirect product decomposition
is the canonical action of the quotient group $\rth\Z$ of $\rth\Z{\times}_{\mu_r}\hs H$ on the normal abelian subgroup $A$.
Since the element $\rth\in\rth\Z$ has image $\zeta$ in $\mu_r$, which lifts to $\varphi^{-1}(\zeta)\in B\subset H$,
we see that $\rth\in\rth\Z$ lifts to $(\rth,\varphi^{-1}(\zeta))\in \rth\Z{\times}_{\mu_r}\hs B\subset \rth\Z{\times}_{\mu_r}\hs H$,
hence $\rth$ acts as $\varphi^{-1}(\zeta)$ on $A$.
Identifying $\rth\Z$ with $\Z$ via $x\mapsto rx$ for $x\in\rth\Z$,   we obtain the  assertion of the proposition.
\qed
\medskip

Now let us take $m=11$, then $A=\Z/11\Z$.
 We take $B=(\Z/11\Z)^{*2}$, the group of nonzero quadratic residues modulo 11.
The group $B$ is a cyclic group of order 5, namely, $B=\{\bar1,\bar4,\bar9,\bar5,\bar3\}$.
Then $H=A\rtimes B$ is a finite nonabelian group of order 55.
Let $n\ge 5$, then  there exists an embedding $\alpha\colon H\into \SL(n,\C)$ .
For $b\in B$, $b\neq \bar1$, let $\varphi_b$ denote the embedding  $B\into\C^*$ taking the generator $b$ of $B$ to $\zeta$, then $\varphi_b^{-1}(\zeta)=b$.
We write $X_{\alpha,b}$ for $X_{\alpha,\varphi_b}$.
Let $\sigma$ be any field automorphism of $\C$ taking $\zeta$ to $\zeta^2$.
Consider the conjugate variety $\sigma X_{\alpha,b}$.

\begin{theorem}\label{t:neq}
For $A=\Z/11\Z$, \ $B=(\Z/11\Z)^{*2}$,   \  $\sigma\in\Aut(\C)$ taking $\zeta$ to $\zeta^2$, the groups
$\pi_1(X_{\alpha,4})$ and $\pi_1(\sigma X_{\alpha,4})$ are not isomorphic.
\end{theorem}

\noindent {\bf Proof.}
We have $\sigma(\zeta)=\zeta^2$.
The homomorphism $\sigma\circ\varphi\colon B\to \C^*$ takes $\bar4$ to $\sigma(\zeta)=\zeta^2$,
hence it takes $\bar4^3=\bar9$ to $(\zeta^2)^3=\zeta$.
Thus $\sigma\circ\varphi_{4}=\varphi_{9}$.

For our group $G$ defined over $\Q$ and for $X=G/\Gamma$, we have
$\sigma X= G/\sigma(\Gamma)$,
where $\sigma$ acts on $\SL(n,\C)$ and on $\C^*$ via the action on $\C$.
For an embedding $\varphi\colon B\into\C^*$  we have
\[\sigma X_{\alpha,\varphi}=G/(\sigma\circ\psi_{\alpha,\varphi})(H)=
G/\psi_{\sigma\circ\alpha,\sigma\circ\varphi}(H)=X_{\sigma\circ\alpha,\sigma\circ\varphi}\,.\]
We obtain that
\[\sigma X_{\alpha,4}=\sigma X_{\alpha,\varphi_{4}}=X_{\sigma\circ\alpha,\sigma\circ\varphi_{4}}=X_{\sigma\circ\alpha,\varphi_{9}}=X_{\sigma\circ\alpha,9}.\]
By Proposition \ref{p:pi-1-X} we have
\[ \pi_1(X_{\alpha,b})\simeq (\Z/11\Z)\rtimes_b \Z,\]
hence
\[\pi_1(X_{\alpha,4})\simeq (\Z/11\Z)\rtimes_{4}\Z\quad\text{and}\quad
\pi_1(\sigma X_{\alpha,4})=\pi_1(X_{\sigma\circ\alpha,9})\simeq (\Z/11\Z)\rtimes_{9}\Z.\]
Now the theorem follows from the next Lemma \ref{l:neq}.
\qed

\begin{lemma}\label{l:neq}
$(\Z/11\Z)\rtimes_{4}\Z\not\simeq (\Z/11\Z)\rtimes_{9}\Z$.
\end{lemma}

We first need the following  group-theoretic fact.

\begin{lemma}\label{p:Gamma-4-9}
Let $A$ be any group without nonzero homomorphisms into $\Z$.
When $\vk$ is an automorphism of $A$, we write $A\rtimes_\vk\Z$ for the semidirect product of $A$ and $\Z$,
where the generator $t=1$ of $\Z$ acts on $A$ by $\vk$.
Fix two automorphisms $\vk_1,\vk_2\in\Aut(A)$, and denote by $\vkbar_i$ the image of $\vk_i$
in the group of outer automorphisms $\Out(A)$.
If $\vkbar_1$ is conjugate to neither $\vkbar_2$ nor $\vkbar_2^{-1}$ in $\Out(A)$, then the semidirect products
$G_1=A\rtimes_{\vk_1} \Z$ and $G_2=A\rtimes_{\vk_2} \Z$ are not isomorphic.
\end{lemma}

\noindent {\bf Proof.}
By contraposition, let  $\lambda\colon G_1\isoto G_2$ be an isomorphism.
Since for each of $i=1,2$, the subgroup $A$ is equal to the kernel
of some/any nonzero homomorphism $G_i\to\Z$, we have $\lambda(A)=A$.
Let $\vk\in\Aut(A)$ denote the restriction of $\lambda$ to $A$.
For the generator $t\in \Z\subset G_1$,
write $\lambda(t)$ as $a\hs t^e\in G_2$ with $a\in A$ and $e\in\Z$.
Since $t$ generates $G_1/A$, we see that $\lambda(t)$ generates $G_2/A$ and hence $e=\pm 1$.
Then for all $a'\in A$, writing $\gamma_a(a')=aa'a^{-1}$ we have
\[\vk(\vk_1(a'))=\lambda(ta't^{-1})=a\hs t^{e} \vk(a')\hs t^{-{e}}a^{-1}=\gamma_a(\vk_2^{e}(\vk(a'))),\]
whence $\vk_1=\vk^{-1}\gamma_a\vk_2^{\hs{e}}\hs\vk$.
Hence $\vkbar_1=\vkbar^{\,-1}\ \vkbar_2^{\,{e}}\ \vkbar$.
Thus $\vkbar_1$ and $\vkbar_2^{\,{e}}$ are conjugate in $\Out(A)$. \qed
\medskip

\noindent {\bf Proof of Lemma \ref{l:neq}.}
We use Lemma \ref{p:Gamma-4-9} when $A=\Z/11\Z$, in which case $\Aut(A)=\Out(A)$,  is abelian and can be identified with $(\Z/11\Z)^*$.
Hence the assumption of Lemma \ref{p:Gamma-4-9} in this case is just that $\vk_1$ and $\vk_2^{\pm 1}$ are distinct as elements of $(\Z/11\Z)^*$.
Here $\vk_1=\bar4$ and $\vk_2=\bar9$.  Since modulo $11$ we have $\bar9\neq \bar4$ and $\bar9^{-1}=\bar5\neq \bar4$,
Lemma \ref{p:Gamma-4-9}  applies and we see that
$(\Z/11\Z)\rtimes_{4}\Z\not\simeq (\Z/11\Z)\rtimes_{9}\Z$.
This completes the proofs of  Lemma \ref{l:neq} and Theorem \ref{t:neq}.
\qed

\bigskip

{\sc Acknowledgements.}
We are grateful to Boris Kunyavski\u\i\ for helpful comments.
This note was completed during a stay of the first-named  author
at the Max-Planck-Institut f\"ur Mathematik, Bonn, and he is grateful to this institute for hospitality, support and
excellent working conditions.

\end{document}